\documentstyle[12pt,leqno,amssymb,amscd,graphics]{amsart}

\newtheorem{thm}{Theorem}[section]
\newtheorem{lemma}[thm]{Lemma}
\newtheorem{prop}[thm]{Proposition}

\theoremstyle{definition}
\newtheorem{dfn}[thm]{Definition}
\theoremstyle{remark}

\newtheorem{rl}[thm]{Rule}
\begin{document}

\newcommand{\ct}{\cite}
\newcommand{\pr}{\protect\ref}
\newcommand{\su}{\subseteq}
\newcommand{\pa}{{\partial}}
\newcommand{\e}{\epsilon}
\newcommand{\es}{{\varnothing}}
\newcommand{\D}{{\mathcal D}}

\newcommand{\A}{{5}}
\newcommand{\AP}{{6}}
\newcommand{\B}{{20}}
\newcommand{\BP}{{21}}
\newcommand{\AB}{{27}}

\newcommand{\C}{{100}}
\newcommand{\CP}{{400}}

\newcounter{numb}

\title{The expected genus of a random chord diagram}

\author{Nathan Linial}
\address{School of Computer Science and Engineering, The Hebrew University of Jerusalem,
Jerusalem 91904, Israel} 
\email{nati@@cs.huji.ac.il}
\urladdr{www.cs.huji.ac.il/$\sim$nati}

\author{Tahl Nowik}
\address{Department of Mathematics, Bar-Ilan University, 
Ramat-Gan 52900, Israel}
\email{tahl@@math.biu.ac.il}
\urladdr{www.math.biu.ac.il/$\sim$tahl}

\date{April 27, 2009}
\begin{abstract}
To any generic curve in an oriented surface there corresponds an oriented chord diagram, 
and any oriented chord diagram may be realized by a curve in some oriented surface. 
The genus of an oriented chord diagram is the minimal genus of an oriented surface 
in which it may be realized. 
Let $g_n$ denote the expected genus of a randomly chosen oriented
chord diagram of order $n$. We show that $g_n$ satisfies:
$$g_n = \frac{n}{2} - \Theta(\ln n).$$

\end{abstract}

\maketitle

\section{Introduction}\label{int}

The study of plane curves dates back to C.F. Gauss. Gauss in \ct{g}
has attached to any plane curve with $n$ double points
a $2n$ letter word, as follows. To each double point attach a letter, and then 
register the letters you encounter as you travel along the curve. 
One obtains a word of length $2n$, where each of the $n$
letters appears precisely twice. Such a word is called a \emph{Gauss word}.
Clearly not any Gauss word may be realized by a plane curve, 
and Gauss has pointed out a necessary condition for it to be realizable.
One can enhance the Gauss word of a curve with a mark on one of the two occurrences of each letter,
signifying that the corresponding strand of the curve crosses the curve at the given double point from right 
to left.
Various necessary and sufficient conditions for a Gauss word to be realizable in the plane
have been eventually given, both in the marked and unmarked settings. 
See \ct{dot},\ct{lom},\ct{rr},\ct{r},\ct{t} and references therein.

To avoid the arbitrariness of the assignment of letters to the different double points,
one can replace Gauss words with chord diagrams. 
To a marked Gauss word corresponds an oriented chord diagram, which is by definition a division of the set
$\{1,\dots,2n\}$ into $n$ ordered pairs. This may be represented as a circle with $2n$ designated points,
and $n$ oriented chords connecting pairs of these points.

Though not every oriented chord diagram may be realized by a curve in the plane, it may be realized 
by a curve in some oriented surface. J.S. Carter in \ct{c} has given a direct construction for the minimal 
genus surface in which a given diagram may be realized, as follows. Take an annulus, which is thought of as 
a regular neighborhood of the curve,
and identify $n$ pairs of regions along the annulus according to the prescription of the diagram.
One obtains an orientable surface $F$ with some $d$ boundary components. 
Capping off the boundary components with discs produces the required minimal genus surface. 
The genus of this surface is 
$g = \frac{1}{2}(n+2-d)$ and so $0 \leq g \leq \frac{1}{2}(n+1)$. 
We will refer to this minimal genus as the genus of the given diagram.
A diagram being realizable in the plane is equivalent to it having genus 0. 

For fixed $n$, we are interested in the distribution of genera of the oriented chord diagrams with $n$ chords, 
and we ask: 
What is the expected genus $g_n$ of a randomly chosen diagram?
We show that the expected genus is very close to the maximal possible genus $\frac{1}{2}(n+1)$, 
in fact, we show that $$g_n = \frac{n}{2} - \Theta(\ln n).$$ 
Since $g = \frac{1}{2}(n+2-d)$,
this is equivalent to showing that the expected number of boundary components is $\Theta(\ln n)$.

Though the above description is geometric, counting the number of boundary components may be 
described in a purely combinatorial manner. Traveling along a boundary component
of our surface corresponds to a walk on the chord diagram according to the following rule.
When moving along the circle of the diagram, and arriving at an end of a chord, 
continue your motion along the chord to its other side. 
If your motion along the chord is in the direction (respectively, against the direction) 
of its orientation, then continue your motion along the circle in the same 
(respectively, opposite) direction as you have moved before entering the chord.
So, to identify a boundary component, one travels along the diagram according
to the above rule, until returning to the starting point. Repeating this process, one obtains all
boundary components.

This combinatorial walk along a chord diagram is reminiscent of the walk 
along the cycles of a permutation. The distribution of the cycles of a random permutation
on $n$ letters is well understood, and the expected number of cycles is also $\Theta(\ln n)$.
Indeed, it is for similar reasons that the expected number of cycles in a permutation and the 
expected number of cycles of the walk along a chord diagram are both $\Theta(\ln n)$, 
though as will be seen, the setting of chord diagrams is substantially more complicated.

Since this problem may be formulated both in a topological and a purely combinatorial manner, 
it may be of interest to both topologists and combinatorialists, and indeed the text is aimed 
for both audiences. 
This work may be viewed as part of the 
recent expanding interest in probabilistic questions in topology, as appears in 
\ct{bhk},\ct{dut},\ct{lim},\ct{mw},\ct{ps}.

\section{Definitions and statement of result}\label{def}

Let $F$ be an oriented surface. A \emph{generic curve} in $F$ is an immersion
$c: S^1 \to F$ for which the only self intersections are transverse double points. 
We fix $n$ once and for all, and call a point of $S^1$ a \emph{dot} 
if it is one of the $2n$ points $\{ e^{\pi i k / n} \ : \ 1 \leq k \leq 2n \}$.
If a generic curve $c$ has $n$ double points, 
then there are $2n$ points in $S^1$ mapped into them,
and we will always assume that these $2n$ points are precisely our $2n$ dots. 
A generic curve with such $n$ double points will be called an $n$-curve.

An oriented chord diagram of order $n$ is a division of the set of dots
into $n$ disjoint ordered pairs. 
One can represent  an oriented chord diagram in the plane, by drawing an
oriented chord connecting each ordered pair, 
where the orientation of the chord represents the order of the pair.
An oriented chord diagram will also be called simply a \emph{diagram}.

We denote the set of all diagrams of order $n$ by $\D_n$, 
and we have $|\D_n| = \frac{(2n)!}{n!}$.
Any $n$-curve determines a diagram $D(c) \in \D_n$ as follows. 
The double points of $c$ 
divide the $2n$ dots into pairs, and the orientation of the surface $F$ induces an 
ordering on each pair, in the following way. If $c(a)=c(b)$ for dots $a,b$, and $c'(b),c'(a)$ is a positive basis
with respect to the orientation of the surface,
then the ordered pair $(a,b)$ is taken. See Figure \pr{f1}. 

\begin{figure}
\scalebox{0.8}{\includegraphics{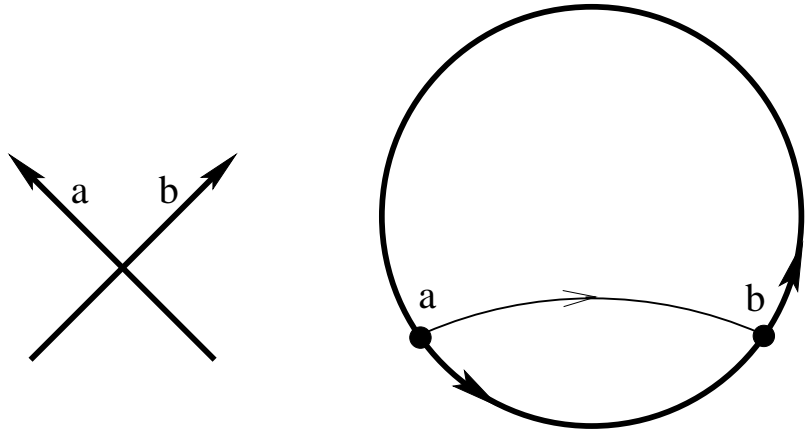}}
\caption{Chord $(a,b)$.}\label{f1}
\end{figure}

Any diagram may be realized by a curve on some 
oriented surface, and a regular neighborhood of the curve in all such surfaces is the same, as we now explain
(compare \ct{c}).
Given a diagram $D \in \D_n$, take the annulus $A = S^1 \times [-\e,\e]$ with a fixed orientation, and identify
$S^1$ with $S^1 \times \{ 0 \}$.
For each dot $a$ let $S_a \su A$ 
be the $2\e \times 2\e$ square centered at $(a,0)$.
We now identify pairs of dots $a,b$ according the prescription of $D$, and we 
identify the corresponding squares $S_a,S_b$ with a positive or negative $\frac{\pi}{2}$ rotation,
so that the self intersection of the curve at the identified square
will be as prescribed by the orientation of the chord between $a$ and $b$ in $D$.
We obtain an oriented surface with $d(D) \geq 1$
boundary components, which we denote $F(D)$. Clearly, 
the embedding $S^1 \to A$ given by $z \mapsto (z,0)$, composed with the quotient
map $A \to F(D)$, is an $n$-curve $S^1 \to F(D)$ which 
realizes $D$, and a regular neighborhood of any curve in any surface realizing $D$, is identical
to $F(D)$. It follows that the genus of the closed surface obtained by capping off the $d(D)$ boundary components
of $F(D)$ with discs, is the minimal genus of a surface in which $D$ may be realized.

\begin{dfn}
Given a diagram $D$, we define $g(D)$, the genus of $D$, to be the minimal genus of
a closed oriented surface admitting a curve $c:S^1 \to F$ with $D(c)=D$.
\end{dfn}

The image in $F(D)$ of our $n$-curve is a graph with $n$ vertices and $2n$ edges,
and so $2-2g(D) = n - 2n + d(D)$, or, $g(D) = \frac{1}{2}(n+2-d(D))$. 
Since $d(D) \geq 1$, we deduce that $0 \leq g(D) \leq \frac{1}{2}(n+1)$.

In this work we study the following question: What is the expected genus of a randomly chosen
diagram $D \in \D_n$? We will show that
the expected genus is close to the maximal possible genus $\frac{1}{2}(n+1)$. 
More precisely we will show:

\begin{thm}\label{main}
The expected genus $g_n$ of a random oriented chord diagram $D \in \D_n$ satisfies:
$$g_n = \frac{n}{2} - \Theta(\ln n).$$
\end{thm}

We think of $g(D)$ and $d(D)$ as random variables defined on a randomly chosen $D \in \D_n$.
That is, our sample space is $\D_n$, each diagram having equal probability $\frac{n!}{(2n)!}$.
We denote the expected values by $g_n = E[g]$ and $d_n = E[d]$.
We will show that $d_n = \Theta(\ln n)$, from which Theorem \pr{main}
follows via $g_n = \frac{1}{2}(n+2-d_n)$.

\section{The random procedure}\label{ran}

We label the dot $e^{\pi i k /n}$, and its corresponding square, 
simply by $k$, and so $k \pm 1$ will mean addition mod $2n$.
The interval along the boundary of the annulus $A$ between two adjacent squares will be 
called an \emph{edge}, so we have $4n$ edges. We orient the edges according to the orientation induced on $\pa A$
from that of $A$, and we denote the oriented edge from square $a$ to square $b$ by $[a,b]$.
The edge $[b,a]$ will then be the parallel edge in the other boundary component. 
So, all edges in $S^1 \times \{ \e \}$ are of the form $[a,a+1]$, 
and they will be called \emph{positive} edges, and all edges in $S^1 \times \{ -\e \}$
are of the form $[a+1,a]$ and will be called \emph{negative} edges.
See Figure \pr{f2}.

\begin{figure}
\scalebox{0.8}{\includegraphics{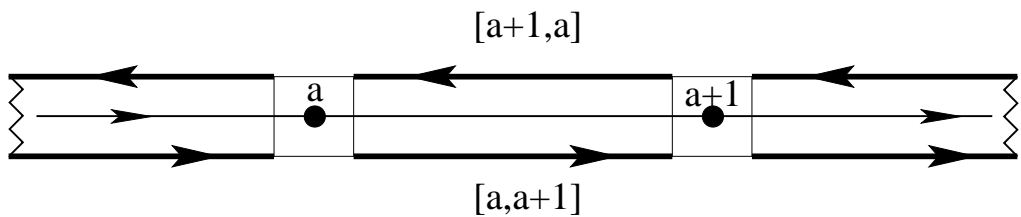}}
\caption{Dots, squares, edges, and their labels.}\label{f2}
\end{figure}

We would like to see how our edges are attached to each other due to the gluing of two squares.
We will say an attachment $[a,b]-[c,d]$ takes place if the end point of $[a,b]$ is glued to the beginning point 
of $[c,d]$. So, say we have identified the two squares $a$ and $b$ according to the oriented chord $(a,b)$.
As seen in Figure \pr{f3}, the eight edges involved are attached to each other as follows:
$[a-1,a]-[b,b+1]$, \ $[b+1,b]-[a,a+1]$, \ $[a+1,a]-[b,b-1]$, \ $[b-1,b]-[a,a-1]$.
This can be summarized by the following rule: 

\begin{figure}
\scalebox{0.8}{\includegraphics{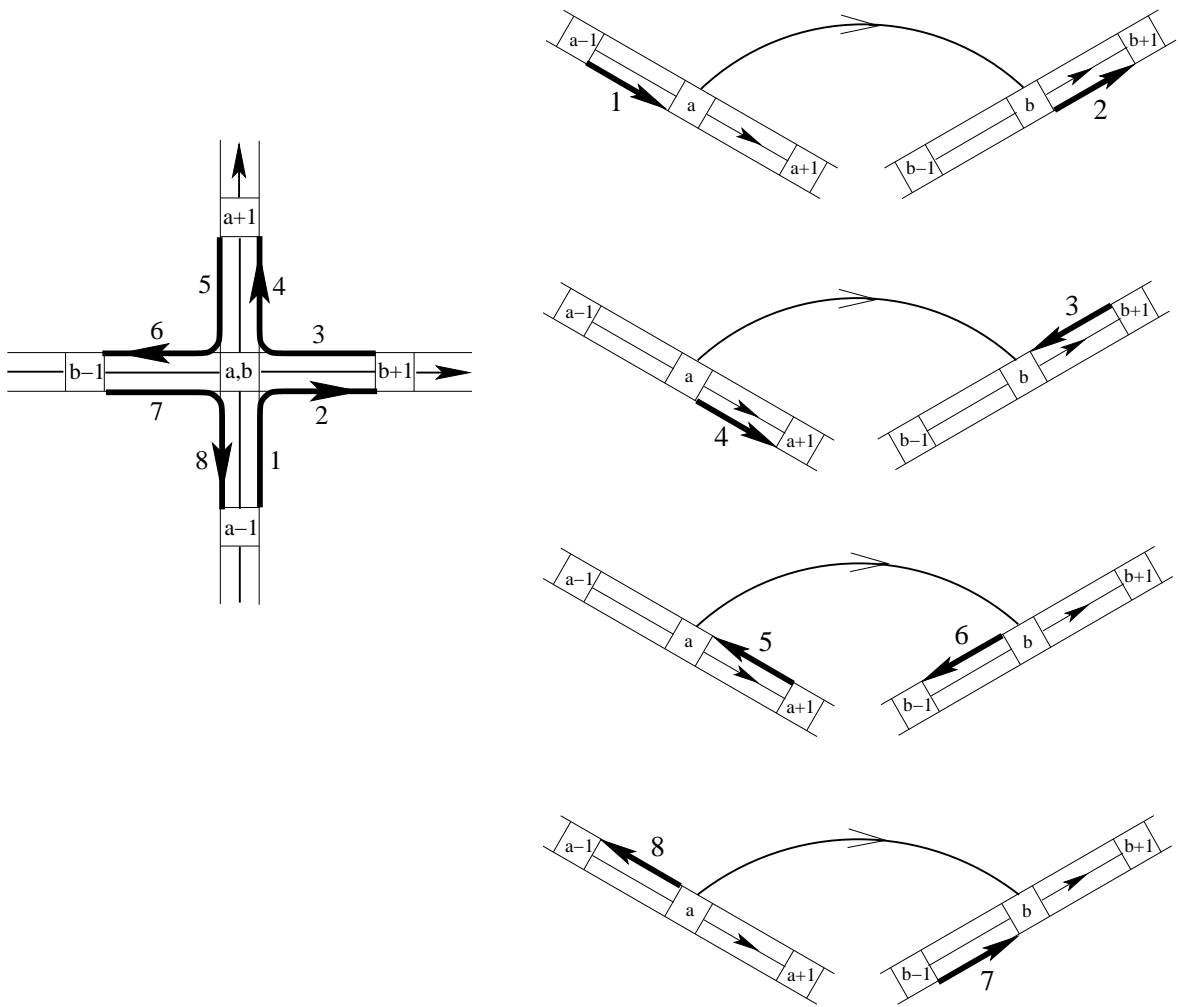}}
\caption{Edge attachments due to chord $(a,b)$.}\label{f3}
\end{figure}

\begin{rl}\label{r}
An attachment $[a,b]-[c,d]$ holds if the oriented chord $(b,c)$ exists, and the signs of $[a,b]$ and $[c,d]$ are the same, 
or if the oriented chord $(c,b)$ exists and the signs of $[a,b]$ and $[c,d]$ are opposite.
\end{rl}

Traveling along a boundary component of the surface $F(D)$ corresponds to
a walk along the diagram $D$, which by Rule \pr{r} proceeds as follows:
When moving along the circle of the diagram, and arriving at an end of a chord, 
continue your motion along the chord to its other side. 
If your motion along the chord is in the direction (respectively, against the direction) 
of its orientation, then continue your motion along the circle in the same 
(respectively, opposite) direction as you have moved before entering the chord.

We may thus read all boundary components directly from the diagram $D$, as follows:
Choose an arbitrary edge, and start traveling along the diagram in the above way, 
alternatingly passing edges and chords, until you return to your initial edge. 
Then choose some unvisited edge, and similarly travel until you return to it.
Continue until all edges have been visited.
Notice that when this is done, each chord of the diagram 
has been visited four times, each visit corresponding to one of the four corners 
of the glued square.

For $k \leq n$, a $k$-$n$-diagram is a choice of $2k$ out of the $2n$ dots, and a division of these $2k$ dots into
$k$ oriented chords, i.e. ordered pairs. The remaining $2n-2k$ dots will be called \emph{vacant} dots. 
So, an oriented chord diagram of order $n$ is an $n$-$n$-diagram.

\begin{dfn}\label{d0}
A \emph{path} in a $k$-$n$-diagram $D$, is a sequence $[a_1,b_1]-[a_2,b_2]- \cdots - [a_r,b_r]$ of distinct edges, 
attached via oriented chords of $D$ according to Rule \pr{r}. 
A path in $D$ is called a \emph{loop} if the attachment $[a_r,b_r]-[a_1,b_1]$ also holds. 
A path in $D$ is called a \emph{segment} if $a_1$ and $b_r$ are vacant dots (perhaps the same dot). 
Loops and segments are precisely those paths that cannot be further extended.
\end{dfn}

We now specify our procedure for choosing a random $n$-$n$-diagram.
Our procedure will choose the chords one by one. 
We first fix an ordering $e_1,\dots,e_{4n}$ of our edges, once and for all. 
Before the procedure begins, we announce $e_1$ as the ``pointer'' edge. Assume that after the $(j-1)$th step, 
we have already chosen $j-1$ oriented chords, and the pointer edge lies in a segment (rather than a loop) 
of the given 
$(j-1)$-$n$-diagram. The next chord is now chosen with one of its dots being the concluding dot of the 
segment in which the pointer lies, 
and its other dot is randomly chosen from the other $2n-2j+1$ vacant dots. The orientation of the new chord
is also randomly chosen. If in the new $j$-$n$-diagram, the pointer's segment continues to be a segment, 
i.e. it does not close into a loop, then the same edge remains the pointer. If on the other hand, 
after the $j$th chord is added, the pointer's segment closes into a loop, then the edge with smallest index
which lies in a segment in the new $j$-$n$-diagram, becomes the new pointer. 
This procedure indeed produces all $n$-$n$-diagrams with equal probability.

Examples of 
two runs of our random procedure appear in Figures \pr{f4} and \pr{f5},
demonstrating some of the interesting features of the 
possible evolution of 
the pointer's segment.
In both figures the edge $e_1$ is the edge $[a,a+1]$ and is marked by 1.

\begin{figure}
\scalebox{0.8}{\includegraphics{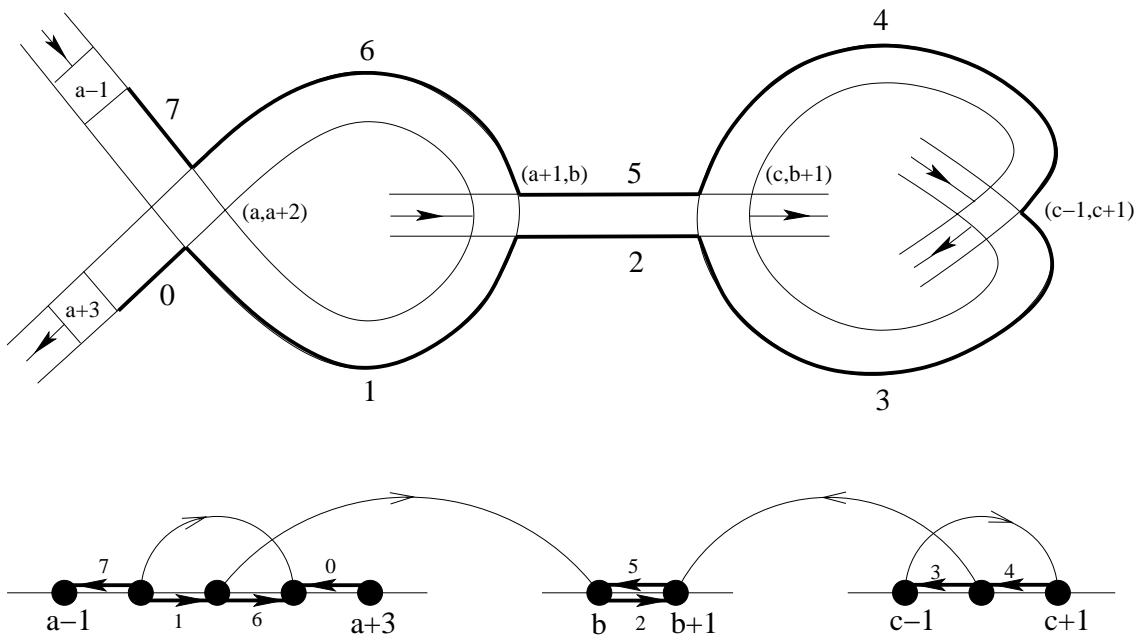}}
\caption{A four step run of the random procedure.}\label{f4}
\end{figure}

\begin{figure}
\scalebox{0.8}{\includegraphics{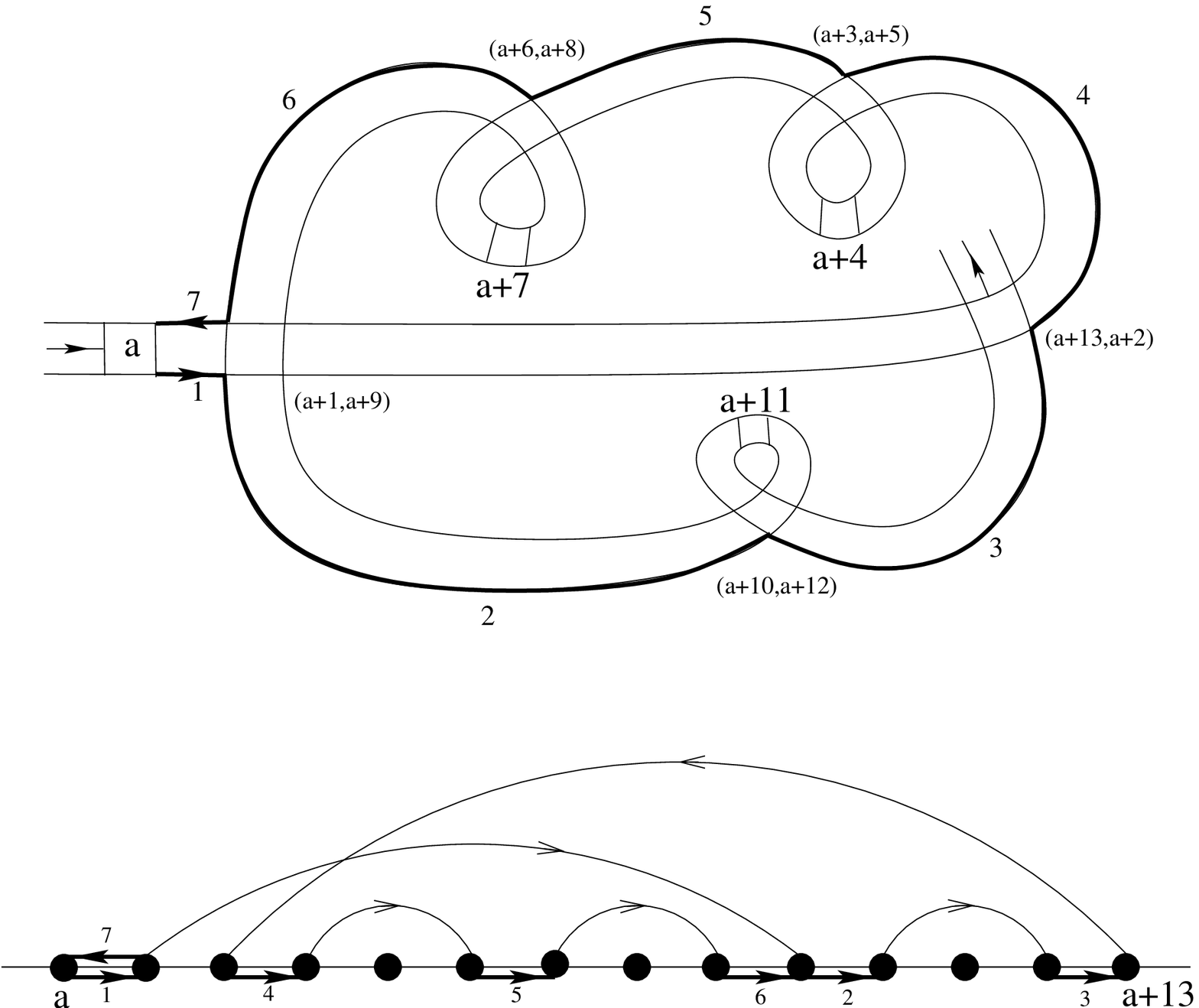}}
\caption{A five step run of the random procedure.}\label{f5}
\end{figure}

In Figure \pr{f4} the chords are chosen in the following order: 
$(a+1,b),(c,b+1),(c-1,c+1),(a,a+2)$.
The segment of $e_1$ 
after the four steps of this run is 
$[a+3,a+2]-[a,a+1]-[b,b+1]-[c,c-1]-[c+1,c]-[b+1,b]-[a+1,a+2]-[a,a-1]$
and these 
edges are marked in the figure by $0,\dots,7$, in this order.
The evolution of $e_1$'s segment throughout the four 
steps of this run
is $1-2$, $1-2-3$,
$1-2-3-4-5-6$, $0-1-2-3-4-5-6-7$.
Note that after the third step the segment 
revisits some chords that 
were chosen in previous steps. Note that after the fourth step the segment 
extends in 
both directions, and so $1$ is no longer the first edge 
in the segment. We point out the following difference between this run, 
and the run of Figure \pr{f5}. In Figure \pr{f4} the initial and final dots of the pointer's segment are 
distinct, whereas in Figure \pr{f5} it is the same dot.

\section{Upper bound for $d_n$}\label{up}

In this section we establish our upper bound for $d_n$. 
If $\ell$ is a loop in an $n$-$n$-diagram $D$, then the \emph{size} of $\ell$ is the number of distinct chords
visited by $\ell$. If $k$ is the size of the loop $\ell$ and $r$ is the number of edges it visits, then since 
$\ell$ alternatingly visits an edge and a chord, and since each chord is visited at most four times,
we have $k \leq r \leq 4k$.

For given $n$, let $L_k = L_k(n)$ denote the expected number of loops of size $k$ in a random $n$-$n$-diagram,
then $d_n = \sum_{k=1}^n L_k$. 
We will show in Proposition \pr{p4} below, that for $k \leq \frac{n}{\C}$, 
$L_k \leq \frac{3}{k}$.  On the other hand, since any chord is visited by at most four different loops,
the total number of all loops of size $k > \frac{n}{\C}$ is at most
$\CP$, and so its expected value $\sum_{k > \frac{n}{\C}} L_k$ is at most $\CP$. Together this gives 
$$d_n = \sum_{k \leq \frac{n}{\C}} L_k + \sum_{k > \frac{n}{\C}} L_k  
\leq \sum_{k \leq \frac{n}{\C}} \frac{3}{k} + \CP \leq 3 \ln n + \CP.$$

In order to obtain our bound $L_k \leq \frac{3}{k}$ we will
need to bound the probability that at a given $k$th step,
the pointer's segment closes into a loop. 
We now prove the following:

\begin{prop}\label{p3}
For $k \leq \frac{n}{\C}$, the probability that the pointer's segment will close 
into a loop at the $k$th step is at most $\frac{3}{4n}$.
\end{prop}

Let $S$ be the pointer's segment in our $(k-1)$-$n$-diagram after step $k-1$, and let
$p$ be its concluding dot. We now need to choose the $k$th chord, with one end being $p$. 
We must determine how many choices will result in closing $S$
into a loop. Let $q$ be the initial dot of $S$. If $q \neq p$ 
(as occurs in the example in Figure \pr{f4}, where $q=a+3$ and $p=a-1$), 
then for $S$ to be closed into a loop, we must choose
$q$ as the second dot for the $k$th chord, and so the choice of unoriented chord is unique.
Though usually only one of the two choices of orientation for this chord will indeed close the
segment $S$ into a loop (as is the case in Figure \pr{f4}), it may in fact occur that both 
orientations accomplish this. Since we are seeking an upper bound for the probability, 
we will always count both orientations as possible, or in other words we will
ignore the choice of orientation in the computation of the probability.
Since there are $2n-2k+1$ vacant dots from which we may choose the second dot for the new chord, the probability 
of choosing the correct dot $q$ is $\frac{1}{2n-2k+1}$, and since we assume $k \leq \frac{n}{100}$ we have 
$\frac{1}{2n-2k+1} \leq \frac{1.1}{2n}$.

If $q=p$ (as in the example of Figure \pr{f5}, where $q=p=a$), then at first sight it may seem that closing the 
pointer's segment into a loop is impossible, and in most cases this is in fact true. But on the other hand, 
there are cases with $q=p$ where not only does there exist a choice of chord that closes the given segment into
a loop, but there are in fact many such choices. The example in Figure \pr{f5} is such case.
The oriented chords $(a,a+4)$, $(a,a+7)$ or $(a,a+11)$ may each be added in the present step to
close the pointer's segment into a loop. For the analysis of this phenomenon, we define the following notion.

\begin{figure}[t]
\scalebox{0.8}{\includegraphics{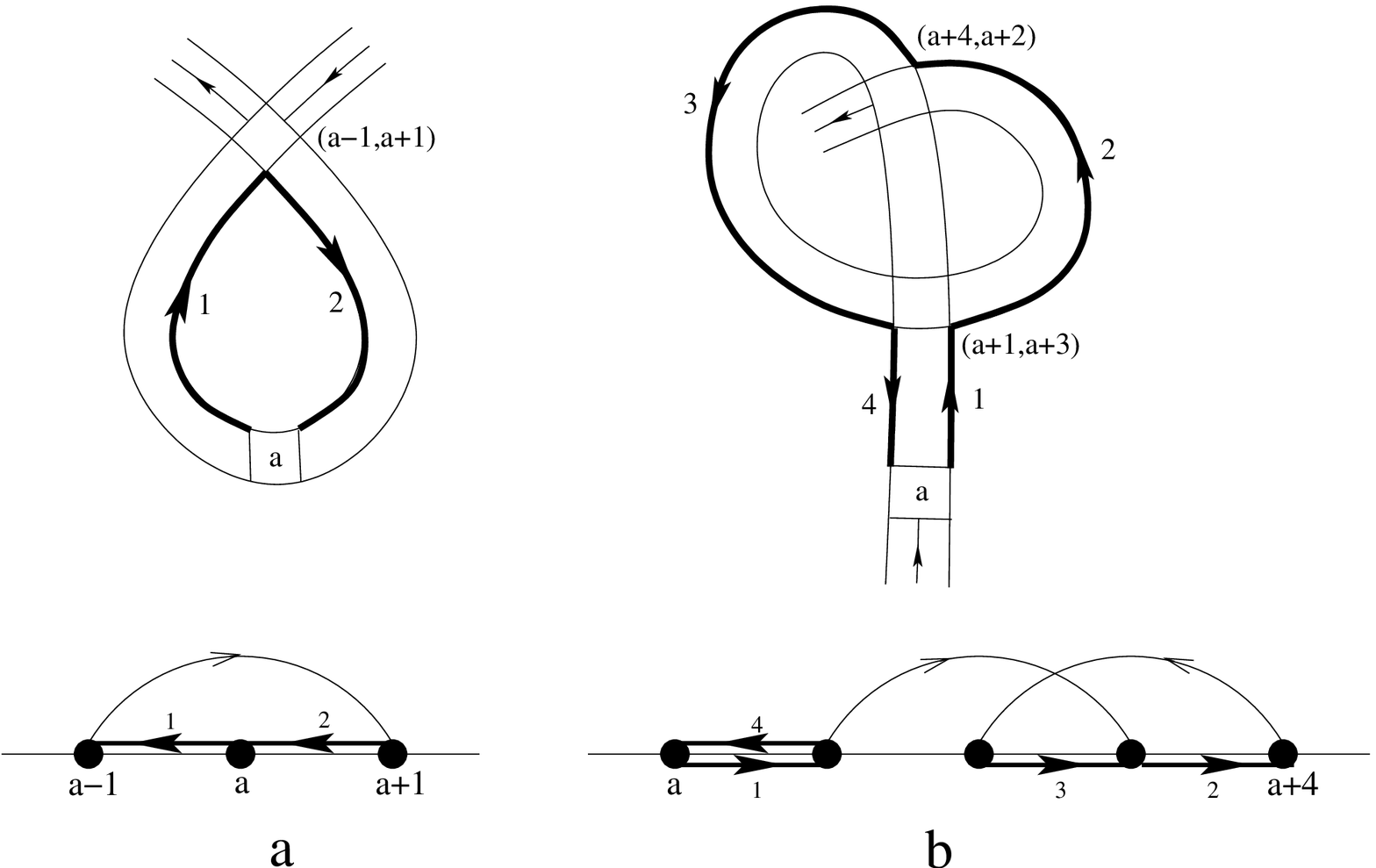}}
\caption{Plugs.}\label{f6}
\end{figure}

\begin{dfn}\label{d1}
A \emph{plug} is a segment $[a_1,b_1]-\cdots-[a_r,b_r]$ for which $a_1=b_r$. 
The vacant dot $a_1=b_r$ is called the \emph{entrance} to the plug.
\end{dfn}

Examples of two plugs are depicted in Figure \pr{f6}.
In Figure \pr{f6}a, the chord $(a-1,a+1)$ produces the plug $[a,a-1] - [a+1,a]$ with entrance $a$. 
In Figure \pr{f6}b, the chords $(a+1,a+3)$ and $(a+4,a+2)$ produce 
the plug $[a,a+1] - [a+3,a+4] - [a+2,a+3] - [a+1,a]$ with entrance $a$. 
Note that the same vacant dot can be the entrance to two different plugs.

In our present case, where $q=p$, the pointer's segment is itself a plug, but this fact is not of interest
to us. What enables us to close the pointer's segment into a loop in Figure \pr{f5}, 
is each one of the \emph{additional} plugs that are present
in the given $5$-$n$-diagram, namely the three plugs $[a+4,a+3]-[a+5,a+4]$, $[a+7,a+6]-[a+8,a+7]$,
and $[a+11,a+10]-[a+12,a+11]$ (each of which is similar to the plug in Figure \pr{f6}a).
In fact, the next lemma shows that in order to close the pointer's segment into a loop in the case $q=p$, 
it is necessary that the second dot of the new chord will be an entrance to a plug.

\begin{lemma}\label{l1}
Let $D$ be a $(j-1)$-$n$-diagram, and let $a$ be a vacant dot in $D$. Let $e$ be an edge entering $a$, 
(i.e. $e$ is $[a-1,a]$ or $[a+1,a]$). Assume $Q$ is an additional oriented chord with one end at $a$ and 
the other end 
at some other vacant dot $b$, such that in the $j$-$n$ diagram obtained by adding $Q$,
the path beginning at $e$ leads to an edge $e'$ which is one of the two exiting edges at $a$ 
(i.e. $e'$ is $[a,a+1]$ or $[a,a-1]$). Then $b$ is the entrance to a plug in $D$.
\end{lemma}

\begin{pf}
Assume $b$ is not an entrance to a plug.
In Figure \pr{f7}, edge $e$ is marked, and the two possibilities for $e'$ are marked.
The path beginning at $e$ passes the point $x$, and in order for it to lead to $e'$, 
it must eventually arrive back into the region depicted in the figure. It does not arrive at $y_1$ or $y_3$ 
since we have assumed that $b$ is not an entrance to a plug. If it arrives at $y_4$ then it closes a loop without
passing either possibility for $e'$. So, it must arrive at $y_2$ as depicted, and it then exits 
our region through point $z$. This time its only possibility for returning is through $y_4$, which as before
prevents it from ever arriving at either possibility for the exiting edge $e'$.
\end{pf}

\begin{figure}[t]
\scalebox{0.8}{\includegraphics{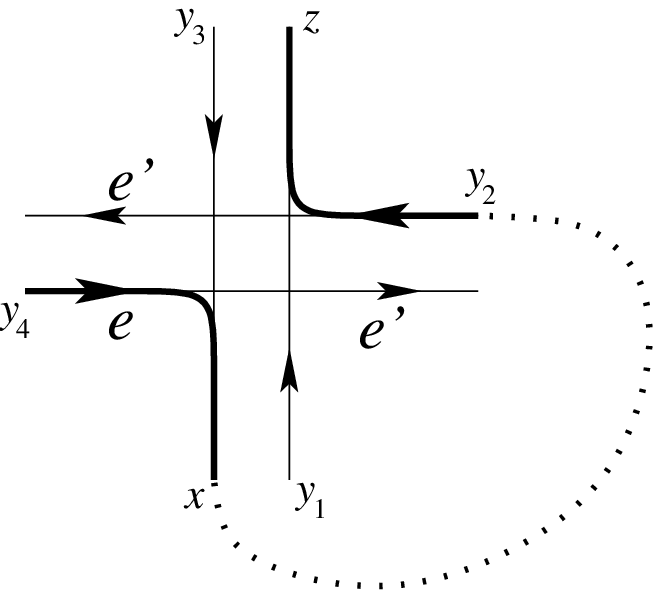}}
\caption{Proof of Lemma \pr{l1}.}\label{f7}
\end{figure}

As we have seen, there may be many plugs available for completing our segment into a loop,
but fortunately, the \emph{expected} number of available plugs is small. The main technical effort of this
work is the following proposition whose proof we defer to Section \pr{plug}.

\begin{prop}\label{pl}
For $k \leq \frac{n}{\C}$, 
the expected number of plugs present after $k$ steps of the random procedure is at most $\frac{1}{4}$. 
\end{prop}

Back to the proof of Proposition \pr{p3} for the case $q=p$. 
By Lemma \pr{l1}, in order for $S$ to close into a loop, the 
second dot we choose for the new chord must be the  entrance to some plug.
By Proposition \pr{pl}, the expected  number of plugs existing 
in the present stage of the random procedure (i.e. after $k-1$ steps), is at most $\frac{1}{4}$.
Together with the case $q \neq p$ we have on average at most $1 + \frac{1}{4}$ choices for the new
unoriented chord. (As before, we ignore the additional choice of orientation.) 
Note that we must take the sum and not the maximum of the
bounds for the two possibilities $q \neq p$ and $q=p$, 
since the expectation for the  number of plugs that we bound
in Proposition \pr{pl} is not conditional on $q = p$. 

We obtain that for $k \leq \frac{n}{100}$, the probability that the
pointer's segment closes into a loop at the $k$th step of the random procedure is at most 
$(1 + \frac{1}{4}) \cdot \frac{1.1}{2n} \leq \frac{3}{4n}$, which completes the proof of Proposition \pr{p3}.

\begin{prop}\label{p4}
Let $L_k$ denote the expected number of loops of size $k$ in a random $n$-$n$-diagram. 
Then for $k \leq \frac{n}{\C}$ we have $L_k \leq \frac{3}{k}$.
\end{prop}

\begin{pf}
It is clear from the definition of our random procedure, that if the segment of $e_1$ closes into a loop
at the $k$th step, then this loop is of size $k$.
So the event that the edge $e_1$ lies in a loop of size $k$ is the same as the event that 
$e_1$ survives as pointer until step $k$, and then at step $k$ its segment closes into a loop.
The probability for this event is at most the probability that at the $k$th step
the pointer's segment closes into a loop, and 
by Proposition \pr{p3} this probability is at most $\frac{3}{4n}$.
Now, our random procedure produces each $n$-$n$-diagram with equal probability, and so by the symmetry of our annulus,
the probability for \emph{any} given edge to lie in a loop of size $k$ is also at most $\frac{3}{4n}$,
or alternatively, the probability $P_k$ that a \emph{randomly} chosen edge will lie
in a loop of size $k$ is at most $\frac{3}{4n}$.
 
We obtain a lower bound for $P_k$ by noting that each loop of size $k$ includes at least $k$ edges, 
and the total number of edges is $4n$, and so $P_k \geq \frac{kL_k}{4n}$. 
Together we get
$\frac{kL_k}{4n}   \leq P_k  \leq \frac{3}{4n}$ which proves our claim.
\end{pf}

As already explained above, the bound $L_k \leq \frac{3}{k}$ for $k \leq \frac{n}{\C}$ 
implies the following upper bound for $d_n$:
$$d_n \leq  3 \ln n + \CP.$$

\section{Lower bound for $d_n$}\label{low}

We have asked in the proof of Proposition \pr{p4},
what is the probability $P_k$ that a randomly chosen edge will lie in a loop of size $k$.
We have noticed that this is precisely the probability that in our random procedure, 
the segment of $e_1$ survives until the $k$th step, 
and then at the $k$th step it closes into a loop. 
In this section we will find a lower bound for $P_k$, for $n \geq 50$ and $k \leq \sqrt{n}$, 
which in turn will provide a lower bound for $d_n$.

A run of $j-1$ steps of the random procedure is called \emph{good} if
after these $j-1$ steps $e_1$ 
is still the pointer, and its segment is 
of the form $[a_1,b_1]-\cdots-[a_j,b_j]$ with all the dots
$a_1,b_1,\dots,a_j,b_j$ 
being distinct. When $j-1=0$, i.e. before starting the random procedure, then the pointer's segment 
is simply $e_1=[a_1,b_1]$, so the run is good.
If the run is good after $j-1$ steps, 
and at the $j$th step the second dot chosen for the new chord is 
not adjacent to any of the dots $a_1,b_1,\dots,a_j,b_j$, then the
run is still good after the $j$th step. This restriction
for the choice of the $j$th chord means that if an edge of the segment is say $[a,a+1]$, then 
the four dots $a-1,a,a+1,a+2$ are not chosen.
So, at the $j$th step we have at most $4j$ dots which we are forbidden to choose,
so the number of allowed choices for a new dot at the $j$th step is at least $2n-4j$. 
Since the total number of dots from which we choose is $2n-2j+1$, 
the probability for such restricted choice at the $j$th step is at least $\frac{2n-4j}{2n-2j+1}$. 

If after $k-1$ steps of the random procedure the run is still good, 
then in particular, the initial and final dots of $e_1$'s 
segment are distinct. So, at the $k$th step there exists 
a choice of oriented chord that closes $e_1$'s segment into a loop, and
the probability for this choice is $\frac{1}{2(2n-2k+1)}$. 
So together, for $n \geq 50$ and $k \leq \sqrt{n}$, the probability that the segment of $e_1$ survives 
until the $k$th step, and then at the $k$th step closes into a loop satisfies:
\begin{gather*}
P_k  \geq    \frac{1}{2(2n-2k+1)}\prod_{j=1}^{k-1} \frac{2n-4j}{2n-2j+1} 
\geq \frac{1}{4n}\prod_{j=1}^{k-1} (1-\frac{2j+1}{2n-2j+1})
\geq \frac{1}{4n}\prod_{j=1}^{k-1} (1-\frac{j+1}{n-k}) \\
\geq       \frac{1}{4n}\prod_{j=1}^{k-1} e^{-\frac{6}{5}\cdot\frac{j+1}{n-k}} 
= \frac{1}{4n}e^{-\frac{6}{5}\sum_{j=1}^{k-1}\frac{j+1}{n-k}}
\geq \frac{1}{4n}e^{-\frac{3}{5}\cdot\frac{k^2 + k}{n-k}} \geq \frac{1}{4n}e^{-\frac{3}{5}\cdot\frac{n + \sqrt{n}}{n-\sqrt{n}}} 
\geq \frac{1}{9n}.
\end{gather*}
(We use the assumption $n \geq 50$ in the fourth and last inequalities.)

As before, let $L_k$ be the expected number of loops of size $k$ in a random $n$-$n$-diagram, then since 
the number of edges in a loop of size $k$ is at most $4k$ we have $P_k \leq \frac{4k L_k}{4n}$.
Together, for $n \geq 50$ and $k\leq \sqrt{n}$ we get
$\frac{1}{9n} \leq P_k \leq \frac{4k L_k}{4n}$, so $L_k \geq \frac{1}{9k}$.
We may now establish our lower bound for $d_n$, for $n \geq 50$:
$$d_n = \sum_{k=1}^{n} L_k \geq \sum_{k=1}^{\sqrt{n}} \frac{1}{9k} \geq \frac{1}{9} \ln \sqrt{n} = \frac{1}{18} \ln n.$$

Together with the upper bound of Section \pr{up} we obtain $d_n = \Theta(\ln n)$, which proves Theorem \pr{main}, 
stating that the expected genus $g_n$ of a random diagram of order $n$ satisfies: 
$$g_n = \frac{n}{2} - \Theta(\ln n).$$

\section{Upper bound for the expected number of plugs}\label{plug}

In this section we prove Proposition \pr{pl}, stating that for $k \leq \frac{n}{100}$, 
the expected number of plugs present in our $k$-$n$-diagram after $k$ steps of the 
random procedure is at most $\frac{1}{4}$.

\begin{dfn}\label{d3}
Two vacant dots in a $k$-$n$-diagram $D$ are called \emph{neighbors}, 
if they are the two end points of a segment in $D$.
\end{dfn}

\begin{dfn}\label{sgn}
A \emph{positive} plug is a plug $[a_1,b_1]-\cdots-[a_r,b_r]$ for which the two edges 
$[a_1,b_1],[a_r,b_r]$ are of the same sign, that is, they are of the form 
$[a,a+1],[a-1,a]$ or $[a,a-1],[a+1,a])$, (as in Figure \pr{f6}a). A \emph{negative}
plug is a plug for which these two edges are of opposite sign, that is, they are of the form 
$[a,a+1],[a+1,a]$ or $[a,a-1],[a-1,a]$, (as in Figure \pr{f6}b). 
Note that if same vacant dot is the entrance to two different plugs, 
then these two plugs must be of the same sign.
\end{dfn}

\begin{lemma}\label{l2}
Under the assumptions of Lemma \pr{l1}, if $e$ and $e'$ are of opposite sign, 
and if $b$ is not the entrance to a positive plug
(and so by Lemma \pr{l1} it is the entrance to one or two negative plugs), 
then either $a$ and $b$ are neighbors (Definition \pr{d3}), or $a$ is also an entrance to a plug.
\end{lemma}

\begin{pf}
Assume  $a$ and $b$ are not neighbors. In order for us to arrive at $e'$,
given that $b$ is not the entrance to a positive plug and $a$ and $b$ are not neighbors,
our path must be as in Figure \pr{f8}, 
which shows that $a$ is the entrance to a (negative) plug. 
\end{pf}

\begin{figure}[t]
\scalebox{0.8}{\includegraphics{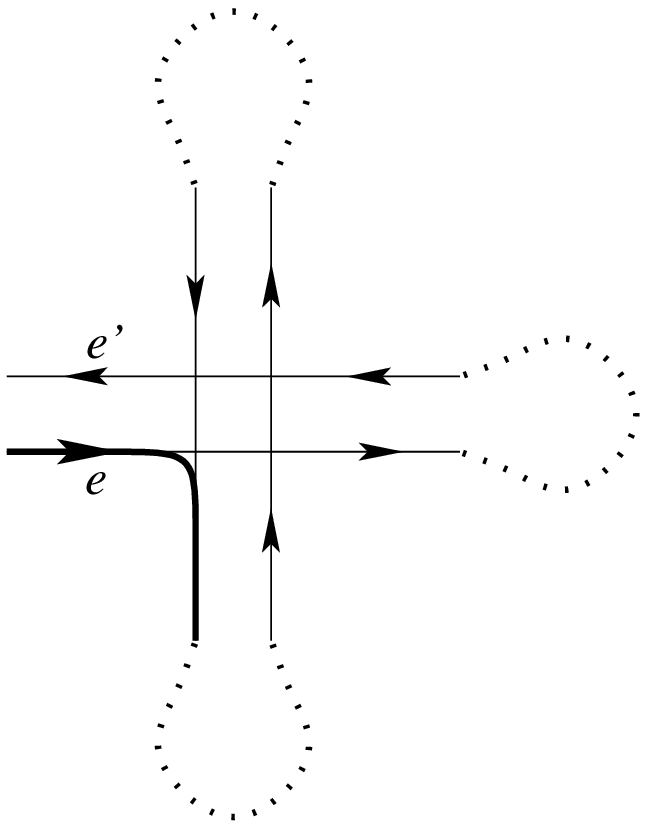}}
\caption{Proof of Lemma \pr{l2}.}\label{f8}
\end{figure}

Any chord is involved in at most four different segments, 
and so at each step, when adding a new chord, at most four new plugs can be created. 
But we will show that in fact the expected number of plugs created at each step $k \leq \frac{n}{\C}$ is at 
most $\frac{25}{n}$. This implies that the expected number of plugs present after $k \leq \frac{n}{\C}$ steps 
is at most $\frac{1}{4}$. 
To establish this bound we will in fact need to prove the following more detailed proposition, 
which distinguishes between positive and negative plugs.

\begin{prop}\label{p2}
The following holds for $k \leq \frac{n}{\C}$:

\begin{enumerate}
\item Let $G^+_k$ (respectively $G^-_k$) denote the expected number of positive (respectively negative) 
plugs completed at the $k$th step.
Then $G^+_k \leq \frac{\A}{n}$ and $G^-_k \leq \frac{\B}{n}$. \label{1}
\item The expected number of plugs present after $k$ steps is at most $\frac{1}{4}$. \label{2}
\item Let $H^+_k$ (respectively $H^-_k$) 
denote the probability that after the $k$th step the concluding dot of the pointer's segment 
is an entrance to a positive (respectively negative) plug.
Then $H^+_k \leq \frac{\AP}{n}$ and $H^-_k \leq \frac{\BP}{n}$. \label{3}
\end{enumerate}
\end{prop}

\begin{pf}
(\pr{1}) Say at the $k$th step we have chosen a chord $Q$ between dots $a$ and $b$, 
and a plug has been completed, with dot $c$ being its entrance. 
This means that after adding $Q$ there is a segment with edges $e_{i_1} - e_{i_2} - \cdots - e_{i_r}$ 
beginning and ending at the vacant dot $c$, and before adding $Q$ this segment did not exist. 
This means that before adding $Q$, 
the segment $S_1$ beginning with $e_{i_1}$ ended at some 
vacant dot $a \neq c$, and the segment $S_2$ ending with $e_{i_r}$ began at some vacant dot
$a' \neq c$. We now distinguish three cases as follows.
If $a \neq a'$ then the new chord $Q$ must be between $a$ and $a'$. 
By definition of our random procedure, 
the concluding vacant dot $p$ of the pointer's segment is one of the dots of the new chord $Q$, 
so must be either $a$ or $a'$.
We will refer to this case as Case A.
If on the other hand $a=a'$ then the new chord $Q$ must be between $a$ and some other vacant dot $b$. 
In this case either $p=a$ or $p=b$, and these two possibilities will be referred to as Case B and Case C, respectively.

For Case A, we note that there are at most four different segments with one end being $p$.
The other end of each such segment is a vacant dot that may be $c$ of the above description. 
For each such $c$ there is a unique second segment with which a configuration $S_1,S_2$ as described above 
may arise for a positive plug, and a unique such second segment for a negative plug. 
Our assumption is that $a \neq a'$ and so for each such configuration there is a unique choice of 
unoriented chord with which such a plug may be created.
As discussed in Section \pr{up}, it may be that both choices of orientation for this chord 
bring to the completion of the plug. So here and in all following cases, we do as we have done
in Section \pr{up}, and include both 
choices in our count by simply ignoring the choice of orientation.
As before, the probability for the correct unoriented chord to be chosen in each case is 
$\frac{1}{2n-2k+1}$ since there are $2n-2k+1$ additional vacant dots, and for $k \leq \frac{n}{100}$
we have $\frac{1}{2n-2k+1} \leq \frac{1.1}{2n}$. 
So, the contribution of this case to $G^+_k$ and $G^-_k$ is
at most  $4 \cdot \frac{1.1}{2n}$. Note that it may be that different configurations in our count are completed
into a plug by the same choice of chord, but by the additivity of expectation, the contributions of
all configurations may be added regardless of the dependence between them. 

In Case B, $Q$ is between the dot $p=a$ and the dot $b$, and by Lemma \pr{l1}, 
$b$ must be an entrance to an existing plug.
We bound all possible contributions that may be from choosing the second dot of the new chord as the entrance 
to an existing plug.
Any new chord may participate in at most four different segments, and so at most four new plugs may be completed. 
By induction, we may use (\pr{2}) of the present proposition for $k-1$
to conclude that on average we have at most $\frac{1}{4}$ existing plugs available. 
So, the contribution is on average at most $4 \cdot \frac{1}{4} \cdot \frac{1.1}{2n} = \frac{1.1}{2n}$.
We cannot determine how this contribution will divide between $G^+_k$ and $G^-_k$ and so we add it to both.

In Case C, $p=b$, and our choice of the second dot $a$ for $Q$ is such that $a$ is part of a configuration
of segments $S_1,S_2$ and dots $c,a$ as described above. 
The segments $S_1$ and $S_2$ may or may not pass chords, 
but there is just one special configuration for $S_1,S_2$ in which both $S_1$ and $S_2$ 
do not pass any chord, namely, the configuration where $a$ and $c$ are adjacent dots
along the annulus, and $S_1,S_2$ are the two edges connecting them. 
If the configuration is not this special one, then necessarily the dot $c$ is adjacent along the annulus 
to a dot that
is the end of one of the $k-1$ existing chords. So there are at most $4(k-1)$ possibilities for such dot.
For each such dot $c$ there are two possibilities for a pair of segments $S_1,S_2$ 
that may give rise to a positive plug,
and two possibilities for a negative plug. Together this gives at most $8k$ possible pairs of segments for 
positive plugs and for negative plugs. Now we note that in order for us to be in Case C,
our dot $p$ must be at the entrance to an existing plug after step $k-1$. By induction we may use (\pr{3}) of 
the present 
proposition for $k-1$ to conclude that the probability for us being in Case C is 
at most $\frac{\AP}{n} + \frac{\BP}{n}$.
And so the contribution of the non-special
configurations to $G^+_k$ and $G^-_k$ is at most 
$8k (\frac{\AP}{n} + \frac{\BP}{n})\frac{1.1}{2n} 
\leq 8 \cdot \frac{n}{\C} \cdot \frac{\AB}{n} \cdot \frac{1.1}{2n}
\leq 3 \cdot \frac{1.1}{2n}$.

For the special configuration, if a plug is completed then it is necessarily a negative plug, so contributes only
to $G^-_k$.
If $p=b$ is the entrance to a positive plug, 
which happens, by induction on part (\pr{3}), with probability at most $\frac{\AP}{n}$,
then we take our bound to be simply the total number of choices $2n-2k+1$ for $a$. 
There may be a special configuration on each side of $a$, and so
the contribution to $G^-_k$ is at most 
$2 \cdot (2n-2k+1) \cdot \frac{\AP}{n} \cdot \frac{1.1}{2n} \leq 24 \cdot \frac{1.1}{2n}$.

If $p=b$ is the entrance to a negative plug, which happens by part (\pr{3}), by induction, with probability at most
$\frac{\BP}{n}$, then by Lemma \pr{l2}, we must choose $a$ which is either a neighbor of $b$ 
or the entrance to a plug.
The dot $b$ has at most 4 neighbors. For each such neighbor $a$ there is 
at most one special configuration that may be
completed into a plug, since it may not be on the side of $a$ where the segment from $b$ arrives at $a$.
So, the contribution of this case is at most 
$4 \cdot \frac{\BP}{n} \cdot \frac{1.1}{2n}$.
The second possibility is that $a$ itself is an entrance to a plug, 
but in Case B above we have already counted all 
possible contributions from connecting $p$ to a dot which is the entrance to an existing plug, 
and so we need not count this again here.
The contribution to $G^-_k$ is thus at most 
$4 \cdot \frac{\BP}{n} \cdot \frac{1.1}{2n} \leq \frac{1.1}{2n}$,  since $n \geq \C$ whenever the assumption
$k \leq \frac{n}{\C}$ is relevant.

We add all contributions for $G^+_k$:
$$G^+_k \leq (4 + 1 + 3) \cdot \frac{1.1}{2n} \leq \frac{\A}{n},$$
and for $G^-_k$:
$$G^-_k \leq 
(4 + 1 + 3 + 24 + 1) \cdot \frac{1.1}{2n}
\leq  \frac{\B}{n}.$$

(\pr{2}) In each step $j \leq k$ on average at most  $\frac{\A}{n} + \frac{\B}{n}$ 
plugs are completed, by (\pr{1}), and so after $k$ steps the expected number of plugs is at most 
$k (\frac{\A}{n} + \frac{\B}{n}) \leq \frac{n}{\C}(\frac{\A}{n} + \frac{\B}{n}) = \frac{1}{4}$.

(\pr{3}) If after the $k$th step, the final dot of our spanning segment is the entrance to a positive plug, 
then this plug may either be one that has existed previously, or one that has just been completed.
If it is a plug that has existed previously, 
then in the $(k-1)$-$n$-diagram
we had before the $k$th step, there is a unique segment $S$ leading to its entrance (which is not the plug itself), 
and let $a$ denote the vacant dot at the beginning of $S$. 
In order for us to end up at the entrance to the given plug after adding the $k$th chord,
this chord must include $a$. 
As before, let $p$ denote the concluding dot of the pointer's segment.
If $p \neq a$ then we have one choice for such unoriented chord. 
If $p=a$ then in order for us to continue into the segment $S$,
then by Lemma \pr{l1} the other dot $b$ of the new chord must be the entrance to an existing plug. 
Together we see that in order for us to land at the entrance of an existing plug, we must choose the second
dot for the new chord either as a dot $a$ as described above,
which is uniquely determined by a plug, or as a dot which is itself the entrance to a plug.
By (\pr{2}) we know that there are on average at most $\frac{1}{4}$ previously existing plugs, and so 
this contributes at most  $2 \cdot \frac{1}{4} \cdot \frac{1.1}{2n}$ to the probability.

On the other hand, the probability that after the $k$th step we have landed at the entrance of a positive plug that 
has just been completed, is at most the probability that such a plug has at all been completed at the $k$th step. 
By (\pr{1}) this probability is at most $\frac{\A}{n}$, 
since the expected number of plugs completed is a bound to the probability that at least one
plug has been completed. 
Together we get $H^+_k \leq  2 \cdot \frac{1}{4} \cdot \frac{1.1}{2n} + \frac{\A}{n} \leq \frac{\AP}{n}$. 
In the same way, using $G^-_k \leq \frac{\B}{n}$ we get $H^-_k \leq \frac{\BP}{n}$
\end{pf}

Recall that what we have actually used from Proposition \pr{p2}
is only part (2), which bounds the total number of plugs. The need for this more detailed analysis is
due to the large contribution of existing positive plugs to the completion of new negative plugs in 
Case C with the special configuration. This required that we separate 
between positive and negative plugs in the inductive proof, with a larger bound for the negative plugs.

\end{document}